\documentclass[12pt]{amsart}

\makeatletter
\@namedef{subjclassname@2020}{\textup{2020} Mathematics Subject Classification}
\makeatother

\newcommand{\tf}[2]{\tfrac{#1}{#2}}
\newcommand{\tff}[2]{\tfrac{#1}{#2}}

\newcommand{\mg}{\infty}

\newcommand{\px}{\pi x}

\newcommand{\dsum}{\di\sum}

\DeclareFontFamily{U}{mathx}{\hyphenchar\font45}
\DeclareFontShape{U}{mathx}{m}{n}{
      <5> <6> <7> <8> <9> <10>
      <10.95> <12> <14.4> <17.28> <20.74> <24.88>
      mathx10
      }{}
\DeclareSymbolFont{mathx}{U}{mathx}{m}{n}
\DeclareFontSubstitution{U}{mathx}{m}{n}
\DeclareMathAccent{\widecheck}{0}{mathx}{"71}

\newcommand{\inrr}{\ensuremath{\in\rr}}

\renewcommand{\kill}[1]{}
\newcommand{\dummy}[1]{\mbox{}}

\makeatletter
\newcommand{\xequal}[2][]{\ext@arrow 0055{\equalfill@}{#1}{#2}}
\def\equalfill@{\arrowfill@\Relbar\Relbar\Relbar}
\makeatother

\newcommand{\xeq}{\xequal}

\renewcommand{\k}{\ensuremath{\ol{\mathrm{P}}}}

\newcommand{\ts}[3]{\left[\,{\di #1}\,\right]^{#2}_{#3}}

\newcommand{\m}{\ensuremath{\infty}}

\renewcommand{\k}[1]{\ensuremath{\left({#1}\right)}}

\newcommand{\bca}{\begin{cases}}
\newcommand{\eca}{\end{cases}}

\newcommand{\mug}{\ensuremath{\infty}}

\newcommand{\logx}{\ensuremath{\log x}}

\newcommand{\ff}[2]{\ensuremath{\di\fr{#1}{#2}}}

\renewcommand{\ss}[3]{\ensuremath{\di\int_{#1}^{#2}{#3}\,dx}}

\newcommand{\bpic}{\begin{picture}}\newcommand{\epic}{\end{picture}}

\newcommand{\beda}{\begin{edaenumerate}}
\newcommand{\eeda}{\end{edaenumerate}}

\newcommand{\cd}{\cdots}

\newcommand{\asx}{\ensuremath{\sin^{-1}x}}

\newcommand{\acx}{\ensuremath{\cos^{-1}x}}
\newcommand{\atx}{\ensuremath{\tan^{-1}x}}

\newcommand{\q}{\quad}

\newcommand{\bq}{\begin{quote}}\newcommand{\eq}{\end{quote}}

\newcommand{\vii}{\\[.2in]}
\newcommand{\rt}{\sqrt}
\newcommand{\be}{\begin{enumerate}}\newcommand{\ee}{\end{enumerate}}
\newcommand{\bce}{\begin{center}}\newcommand{\ece}{\end{center}}
\newcommand{\bde}{\begin{description}}\newcommand{\ede}{\end{description}}
\newcommand{\bri}{\begin{flushright}}\newcommand{\eri}{\end{flushright}}
\newcommand{\bb}{\begin{block}}\newcommand{\eb}{\end{block}}
\newcommand{\bt}{\begin{thm}}\newcommand{\et}{\end{thm}}
\newcommand{\bpf}{\begin{proof}}\newcommand{\epf}{\end{proof}}
\newcommand{\bex}{\begin{ex}}\newcommand{\eex}{\end{ex}}
\newcommand{\bexr}{\begin{exr}}\newcommand{\eexr}{\end{exr}}
\newcommand{\bft}{\begin{fact}}\newcommand{\eft}{\end{fact}}
\newcommand{\brk}{\begin{rmk}}\newcommand{\erk}{\end{rmk}}
\newcommand{\ba}{\begin{align*}}\newcommand{\ea}{\end{align*}}
\newcommand{\bexe}{\begin{exe}}\newcommand{\eexe}{\end{exe}}

\newcommand{\bit}{\begin{itemize}}\newcommand{\eit}{\end{itemize}}

\newcommand{\bcm}{}
\newcommand{\ol}{\overline}
\newcommand{\hf}{\hfill}
\newcommand{\fr}{\frac}
\newcommand{\cc}{\ensuremath{\mathbf{C}}}

\newcommand{\qq}{\ensuremath{\mathbf{Q}}}
\newcommand{\rr}{\ensuremath{\mathbf{R}}}

\newcommand{\bd}{\begin{defn}}\newcommand{\ed}{\end{defn}}
\newcommand{\bp}{\begin{prop}}\newcommand{\ep}{\end{prop}}
\newcommand{\p}{\ensuremath{\pi}}
\newcommand{\eh}{\emph}

\newcommand{\te}{\text}

\newcommand{\di}{\displaystyle}

\newcommand{\f}{\frac}

\newcommand{\z}{\ensuremath{\bm{z}}}

\newcommand{\kref}[1]{(\ref{#1})}
\usepackage[dvipdfmx]{graphicx}
\usepackage[dvipsnames]{xcolor}
\usepackage{float,afterpage}
\usepackage{asymptote,layout}
\usepackage{wrapfig,epic}

\usepackage{geometry}
\usepackage{exscale,latexsym,bm}
\usepackage{amssymb,enumerate,amsmath,amsthm,amsfonts}
\usepackage{verbatim,fancybox,wasysym,fancyhdr,type1cm}
\usepackage[frame,all,poly,curve,knot,arrow]{xy}
\usepackage{colortbl}
\usepackage{boxedminipage}
\usepackage{multirow}

\graphicspath{{./figs/}}
\theoremstyle{definition}
\newtheorem{thm}{Theorem}[section]

\newtheorem{lem}[thm]{Lemma}
\newtheorem{prop}[thm]{Proposition}\newtheorem{cor}[thm]{Corollary}
\newtheorem{cj}[thm]{Conjecture}

\newtheorem{exr}[thm]{Exercise}

\newtheorem{ex}[thm]{Example}

\newtheorem{defn}[thm]{Definition}\newtheorem{rmk}[thm]{Remark}
\newtheorem{fact}[thm]{Fact}
\newtheorem{block}[thm]{}
\newtheorem*{exe}{Exercise}

\renewcommand{\z}{\zeta}

\renewcommand{\asx}{\sin^{-1}x}
\renewcommand{\acx}{\cos^{-1}x}

\renewcommand{\kref}[1]{(\ref{#1})}
\renewcommand{\asx}{\sin^{-1}x}
\renewcommand{\acx}{\cos^{-1}x}
\renewcommand{\z}{\zeta}

\renewcommand{\asx}{\arcsin x}
\renewcommand{\acx}{\arccos x}

\renewcommand{\atx}{\arctan x}

\newcommand{\athx}{\text{arctanh}\,x}
\newcommand{\arcsinh}{\text{arcsinh}\,}
\newcommand{\arctanh}{\text{arctanh}\,}

\newcommand{\logrx}{\log^{r}x}
\newcommand{\cbn}{\ff{\binom{2n}{n}}{2^{2n}}}
\newcommand{\cbk}{\ff{\binom{2k}{k}}{2^{2k}}}
\newcommand{\arccot}{\text{\,arccot\,}}
\newcommand{\arctanx}{\atx}
\newcommand{\pix}{\pi x}

\begin{document}
\title[From Basel Problem to multiple zeta values]{
From Basel Problem to multiple zeta values}


\author{Masato Kobayashi}

\address{Masato Kobayashi\\
Department of Engineering\\
Kanagawa University, 3-27-1 Rokkaku-bashi, Yokohama 221-8686, Japan.
\\
}
\email{masato210@gmail.com}




\date{December, 5, 2021.}                             




\begin{abstract}
We show many consequences of the proof of Basel problem by Boo Rim Choe (1987) to central binomial series, multiple zeta values and some other multiple sums. The main idea is to evaluate integrals involving powers of arcsine function.
\end{abstract}

\maketitle
\tableofcontents


 

\section{Introduction}
%

\subsection{Basel Problem: arcsin and $\z(2)$}

Let $\z$ denote the Riemann zeta function.
Boo Rim Choe (1987) \cite{boo} gave evaluation of the integral 
\begin{eqnarray}
\ff{\p^{2}}{8}=
\di\int_{0}^{1}{\ff{\arcsin x}{\rt{1-x^{2}}}}\,dx	=
\ff{3}{4}\z(2)
\label{ba1}
\end{eqnarray}
which provides another proof of \emph{Basel problem} $\z(2)=\tfrac{\p^{2}}{6}$ (dating back to Euler around 1735).
Actually, there is a counterpart of this:
\begin{eqnarray}
\ff{2}{\pi} 
\di\int_{0}^{1}{\ff{\arcsin^{2} x}{2!}}\,
\ff{dx}{\rt{1-x^{2}}}
=\ff{1}{4}\z(2)
\label{ba2}
\end{eqnarray}
as the even part of $\z(2)$;
we will explain why this equation involes $\tfrac{2}{\p}$ and $\arcsin^{2}x/2!$ later on. The aim of this article is to show 
analogous integral evaluation 
\begin{eqnarray}
\di\int_{0}^{1}\ff{\asx\acx}{x}dx=\f{7}{8}\z(3)
\label{eq1}
\end{eqnarray}
\begin{eqnarray}
\ff{2}{\pi}\di\int_{0}^{1}
\ff{\arcsin^{2}x}{2!}
\ff{\acx}{x}dx
=\f{1}{8}\z(3)
\label{eq2}
\end{eqnarray}
and discuss its applications 
to central binomial series, multiple zeta values and some other sums shown as 
Theorems \ref{thm1}, \ref{thm2}, \ref{thm3}, \ref{thm4}, \ref{thm5}, \ref{thm6}.

\begin{rmk}
As the referee kindly pointed out, Mathematica actually returns \kref{eq1}.
However, we could not find any reference in this context which particularly deals with appearance of arccos. 
In this sense, our result is new.
\end{rmk}




\subsection{central binomial sum, Wallis integral}

One of important topics in number theory 
is \emph{central binomial sums}.
Informally speaking, it is  an infinite series involving  $\binom{2n}{n}$.
Lehmer \cite{le} discussed two types of such sums
\[\textrm{I}.\, 
\sum_{n=0}^{+\infty} a_n \binom{2n}{n}, \quad 
\textrm{II}.\,
\sum_{n=0}^{+\infty} a_n \binom{2n}{n}^{-1}.
\]
Some examples are 
\[
\dsum_{n=0}^{\mug}\ff{\binom{2n}{n}}{8^{n}}
=\rt{2}, \q 
\dsum_{n=0}^{\mug}\ff{\binom{2n}{n}}{10^{n}}
=\rt{\ff{5}{3}}.
\]
He also presented connection between such series and 
Maclaurin series of $\arcsin x$ and $\arcsin^2x$.
Other examples are 
\[
	\sum_{n=1}^{\mg}\ff{1}{n^2\binom{2n}{n}}=\ff{1}{3}\,\z(2), \quad
	\sum_{n=1}^{\mg}\ff{(-1)^{n-1}}{n^3\binom{2n}{n}}=\ff{2}{5}\,\z(3),
\]
as they arise in the work of Ap\`{e}ry \cite{apery} and van der Poorten \cite{van} to prove irrationality of $\z(3)$. 


\begin{rmk}
For a nonnegative integer $n$, let
\begin{align*}
	(2n-1)!!&=(2n-1)(2n-3)\cd 3\cdot 1,
	\\(2n)!!&=2n(2n-2)\cd 4\cdot 2.
\end{align*}
We understand $(-1)!!=0!!=1$. Notice that the following relation holds.
\[
\cbn= \ff{(2n-1)!!}{(2n)!!}.
\]
These numbers appear in Wallis integral:
\[
\ss{0}{\pi/2}{\sin^n x}
=
\begin{cases}
	\ff{\p}{2} \cbn&	\text{$n$ even,}\vii
	\cbn&	\text{$n$ odd.}\\
\end{cases}
\]
\end{rmk}


\section{main results}

\subsection{arcsin and $\z(3)$}

Toward the proof of \kref{eq1}, \kref{eq2}, 
we first setup some notation for convenience.


\begin{defn}
Let $\rr[[x]]$ denote the set of real power series. For $f(x)\in \rr[[x]]$, define $W:\rr[[x]]\to \rr[[x]]$ by 
\[
Wf(x)=
\di\int_{0}^{1}{f(xu)\ff{du}{\rt{1-u^{2}}}}.
\]
\end{defn}
In particular, 
\[
Wf(x)
\bigr|_{x=1}=
\di\int_{0}^{1}{f(u)\ff{du}{\rt{1-u^{2}}}}.
\]

\begin{fact}
Recall from calculus that 
if $f(x)=
\sum_{n=0}^{\mug}a_{n}x^{n}$ ($a_{n}\in\rr$) is a convergent power series with the radius of convergence $R$, then so is $\di\int_{0}^{x}{f(u)}du$ and moreover it is given by termwise integration
\[
\di\int_{0}^{x}{f(u)}\,du=
\sum_{n=0}^{\mug}\ff{a_{n}}{n+1}x^{n+1}.
\]
In the sequel, we will use this result without mentioning explicitly.
\end{fact}


\begin{lem}
\label{le1}
Let $f\left(x\right)\in \mathbf{R}[[x]]$. 
\begin{enumerate}
	\item Moreover, suppose it is odd in the form 
\[f\left(x\right)= \sum_{k=0}^{\mg}
\cbk
\,a_{2k+1}\,x^{2k+1}, \q 
a_{2k+1}\inrr.
\]
Then 
\[Wf\left(x\right)=\sum_{k=0}^{\mg}\f{a_{2k+1}}{2k+1}\,x^{2k+1}.\]
	\item Moreover, suppose it is even in the form 
\[f\left(x\right)= \sum_{k=0}^{\mg}
{
\ff{2^{2k}}{\binom{2k}{k}}
}
\,a_{2k}\,x^{2k}, \q a_{2k}\inrr.
\]
Then 
\[Wf\left(x\right)=
\f{\,\pi\,}{2}
\sum_{k=0}^{\mg}a_{2k}\,x^{2k}.\]
\end{enumerate}

\end{lem}
Thus, by the opertator $W$ we can ``kill" $\cbk$ 
or $\ff{2^{2k}}{\binom{2k}{k}}$ 
from coefficients and instead $\tfrac{1}{2k+1}$ 
or $\tfrac{\pi}{2}$ shows up, respectively.

\begin{proof}
To show (1), recall that 
\[
\di\int_{0}^{1}{\ff{u^{2k+1}}{\rt{1-u^{2}}}}\,du
=\ff{(2k)!!}{(2k+1)!!}
=
\ff{2^{2k}}{\binom{2k}{k}}\ff{1}{2k+1}
.
\]
Then 
\[
Wf(x)=
\int_0^1 \sum_{k=0}^\infty \cbk
a_{2k+1}x^{2k+1}\f{u^{2k+1}}{\sqrt{1-u^2}} du
\]
\[
=
\sum_{k=0}^\infty
\cbk
a_{2k+1}x^{2k+1}
{\int_0^1
\f{u^{2k+1}}{\sqrt{1-u^2}} du}=
\sum_{k=0}^{\mg}
\f{a_{2k+1}}{2k+1}\,x^{2k+1}.\]
We can verify (2) likewise with 
\[
\di\int_{0}^{1}{\ff{u^{2k}}{\rt{1-u^{2}}}}\,du
=\ff{\pi}{2}\ff{\binom{2k}{k}}{2^{2k}}
.
\]
\end{proof}

\begin{lem}\label{le2}
Let $f(x)\in \rr[[x]]$. Suppose moreover $f(0)=0$. Then 
\begin{eqnarray}
W\k{
\di\int_{0}^{x}{\ff{f(y)}{y}\,dy}}
=\di\int_{0}^{1}
	{\ff{f(xu)}{u}}\arccos{u} \,du.
\label{eqle1}
\end{eqnarray}
In particular, for $x=1$, we have 
\[
W\k{
\di\int_{0}^{x}{\ff{f(y)}{y}}\,dy}
\Bigr|_{x=1}
=\di\int_{0}^{1}{\ff{f(u)}{u}}\arccos{u} \,du.
\]
\end{lem}
\begin{proof}
If $x=0$, then both sides in \kref{eqle1} are 0. 
Suppose $x\ne0$. Exchanging order of the double integral (Fubini's Theorem),  we have 
\begin{align*}
	W\k{
\di\int_{0}^{x}{\ff{f(y)}{y}}\,dy}
&=
\di\int_{0}^{1}{
\di\int_{0}^{xu}{\ff{f(y)}{y}}\,dy
}\,\ff{du}{\rt{1-u^{2}}}
	\\&=
\di\int_{0}^{x}{
\di\int_{y/x}^{1}{\ff{f(y)}{y}}
\ff{1}{\rt{1-u^{2}}}
}\,dudy
	\\&=\di\int_{0}^{x}
	{\ff{f(y)}{y}}\arccos{\ff{y}{x}} \,dy
	\\&=\di\int_{0}^{1}
	{\ff{f(xu)}{u}}\arccos{u} \,du.
\end{align*}
\end{proof}

\begin{thm}\label{thm1}
\begin{eqnarray}
\di\int_{0}^{1}\ff{\asx\acx}{x}dx
=\ff{7}{8}\z(3),
\label{eq3}
\end{eqnarray}
\begin{eqnarray}
\ff{2}{\pi}
\int_{0}^{1}
\frac{\arcsin x}{2!}
\ff{\arccos x}{x} dx
=\ff{1}{8}\z(3)
.
\label{eq4}
\end{eqnarray}
\end{thm}

\begin{proof}
Start with the well-known Maclaurin series
\[
\arcsin y=
\dsum_{k=0}^{\mug}\cbk
\ff{y^{2k+1}}{2k+1}, \q |y|\le 1.
\label{1}
\]
Then 
\[
\di\int_{0}^{x}{
\ff{\arcsin y}{y}
}\,dy
=
\dsum_{k=0}^{\mug}\cbk\ff{x^{2k+1}}{(2k+1)^{2}}.
\]
Notice that $f(y)=\arcsin y$ is odd. 
By Lemma \ref{le1} (1),
\[
W
\k{\di\int_{0}^{x}{
\ff{\arcsin y}{y}
}\,dy}
=
\dsum_{k=0}^{\mug}\ff{x^{2k+1}}{(2k+1)^{3}}.
\]
Then $x=1$ yields 
\[
W
\k{\di\int_{0}^{x}{
\ff{\arcsin y}{y}
}\,dy}
\Bigr|_{x=1}
=
\dsum_{k=0}^{\mug}\ff{1}{(2k+1)^{3}}
=\frac{7}{8}\z(3).
\]
By Lemma \ref{le2} and $f(0)=0$, we have 
\[
\di\int_{0}^{1}\ff{\arcsin u\arccos u}{u}du=\frac{7}{8}\z(3).
\]
In the same way, we can prove \kref{eq4} with the 
less-known Maclaurin series 
\[
\ff{\arcsin^{2}x}{2!}=
\dsum_{k=0}^{\mug}\ff{2^{2k}}{\binom{2k}{k}}\ff{x^{2k}}{(2k)^{2}}.
\]
\end{proof}



\subsection{singed analog}
From now on, 
we discuss many consequences of Theorem \ref{thm1}. 
The first is signed analog of \kref{ba1} and \kref{eq1}
with replacing 
$\arcsin x$ by $\arcsinh x$. 


\begin{thm}\label{thm2}
Let $G$ denote the Catalan constant 
\[
1-\ff{1}{3^{2}}+\ff{1}{5^{2}}-\ff{1}{7^{2}}+\cd=
0.915966\cd.
\]
Then 
\[
\di\int_{0}^{1}{\ff{\arcsinh x}{\rt{1-x^{2}}}}\,dx	=
G,
\]
\[
\di\int_{0}^{1}\ff{\arcsinh x\acx}{x}dx=
\ff{\p^{3}}{32}.
\]
\end{thm}
\begin{proof}
Recall from complex analysis that 
\[
\arcsin y=
\dsum_{k=0}^{\mug}
\cbk\ff{y^{2k+1}}{2k+1}, \q |y|\le 1
\]
and $\arcsinh z=-i\arcsin(iz)$ (for all $z\in\cc$). 
Then it follows that 
\[
	\arcsinh y=
	\dsum_{k=0}^{\mug}(-1)^{k}
	\cbk\ff{y^{2k+1}}{2k+1}
\]
for all $y\inrr$ such that $|y|\le 1$ so that 
\[
\di\int_{0}^{1}{\ff{\arcsinh u}{\rt{1-u^{2}}}}\,du
=W\k{\arcsinh x} \Bigr|_{x=1}
=\dsum_{k=0}^{\mug}\ff{(-1)^{k}}{(2k+1)^{2}}=G.
\]
Next, consider the argument with replacing $\arcsin y$ by $\arcsinh y$ in Proof of Theorem \ref{thm1} 
throughout. Then, the proof goes without any substantial changes.
 It leads us to 
%
\[
\di\int_{0}^{1}\ff{\arcsinh x\acx}{x}dx=
\dsum_{k=0}^{\mug}\ff{(-1)^{k}}{(2k+1)^{3}}.
\]
The sums on right hand side is $\f{\p^{3}}{32}$. 
\end{proof}


\subsection{arctan}

The next idea is to replace 
$\arcsin, \arccos$ in \kref{eq3} by 
$\arctan, \text{arccot}$ and see what happens. 
This idea might sound simplistic.
However, the integral happens to remain the same value $\tf{7}{8}\z(3)$. Furthermore, its proof suggests several applications to multiple sums 
in the sequel.

\begin{thm}\label{thm3}
\[
\di\int_{0}^{1}{\ff{\arctan x\arccot x}{x}}\,dx
=\ff{7}{8}\z(3).
\]
\end{thm}
\begin{proof}
Let 
\[
A=
\di\int_{0}^{1}{\ff{\arctanx\arccot x}{x}}\,dx,
\]
\[
A(1)=
\di\int_{0}^{1}{\ff{\arctanx}{x}}\,dx,
\q 
A(2)=
\di\int_{0}^{1}{\ff{\arctan^{2}x}{x}}\,dx.
\]
Then 
\begin{align*}
	A&=\di\int_{0}^{1}{\ff{\arctanx\arccot x}{x}}\,dx
	\\&=\di\int_{0}^{1}{\ff{\arctanx\k{
\f{\p}{2}-\arctanx}
}{x}}\,dx
=
{\ff{\p}{2}A(1)-A(2)}.
\end{align*}
We evaluate $A(1)$ and $A(2)$ as follows.
\begin{align*}
	A(1)&=\di\int_{0}^{1}{\ff{\arctanx}{x}}\,dx
	=\di\int_{0}^{1}{
\dsum_{k=0}^{\mug}\ff{(-1)^{k}}{2k+1}x^{2k}
}\,dx
	\\&=\dsum_{k=0}^{\mug}\ff{(-1)^{k}}{2k+1}
\di\int_{0}^{1}{x^{2k}}\,dx=
\dsum_{k=0}^{\mug}\ff{(-1)^{k}}{(2k+1)^{2}}=
G.
\end{align*}
For $A(2)$, recall from Fourier analysis that 
\[
\log\k{\tan\ff{y}{2}}=-2
\dsum_{k=0}^{\mug}\ff{1}{2k+1}\cos(2k+1)y, \q 0<y<\pi.
\]
It follows that 
\begin{align*}
	A(2)&=\di\int_{0}^{1}{\ff{\arctan^{2}x}{x}}\,dx
\xeq[y=2\atx]{}
\ff{\,1\,}{4}
\di\int_{0}^{\p/2}{\ff{y^{2}}{\sin y}}\,dy
	\\&=\ff{\,1\,}{4}
\k{\ts{y^{2}
\log\k{\tan \ff{y}{2}}
}{\p/2}{0}
-
\di\int_{0}^{\p/2}{
2y\log\k{\tan \ff{y}{2}}
}\,dy
}
	\\&=-\ff{1}{2}
\di\int_{0}^{\p/2}{
y
\k{-2
\dsum_{k=0}^{\mug}\ff{1}{2k+1}
\cos(2k+1)y
}
}\,dy
	\\&=\dsum_{k=0}^{\mug}
\ff{1}{2k+1}
\di\int_{0}^{\p/2}{y\cos(2k+1)y}\,dy
	\\&=\dsum_{k=0}^{\mug}
\ff{1}{2k+1}
\k{\ts{
y\ff{\sin(2k+1)y}{2k+1}
}{\p/2}{0}-
\di\int_{0}^{\p/2}{
\ff{\sin(2k+1)y}{2k+1}
}\,dy
}
	\\&=\dsum_{k=0}^{\mug}
\k{\ff{\,\pi\,}{2}\ff{(-1)^{k}}{(2k+1)^{2}}
-\ff{1}{(2k+1)^{3}}}
=\ff{\pi G}{2}-\ff{7}{8}\z(3).
\end{align*}
Conclude that 
\[
A=
{
\ff{\p G}{2}-
\k{\ff{\pi G}{2}-\ff{7}{8}\z(3)
}
}=\ff{7}{8}\z(3).
\]
\end{proof}

\subsection{multiple values}

As a natural generalization of Riemann zeta function,
let us introduce the following sums.

\begin{defn}
For positive integers $i_1, \dots, i_k$ such that $i_1\geqslant 2$, define the \emph{multiple zeta value} and 
\emph{multiple $t$-value} by 
\[
\z\left(i_1, i_2, \ldots, i_k\right)
=
\sum_{
\substack{
n_1>n_2>\cdots >n_k
}
}
\frac{1}{n_1^{i_1}n_2^{i_2}\cdots n_k^{i_k}},
\]
\[
t\left(i_1, i_2, \ldots, i_k\right)
=
\sum_{
\substack{
n_1>n_2>\cdots >n_k\\
n_j\, \mathrm{\,odd}
}
}
\frac{1}{n_1^{i_1}n_2^{i_2}\cdots n_k^{i_k}}.
\]
\end{defn}

Sometimes it is better to 
interpret $\z(i_{1}, \cd, i_{k})$ as 
\[
2^{i_{1}+\cd+i_{k}}
\sum_{
\substack{
m_1>m_2>\cdots >m_k\\
m_j\, \mathrm{\,even}
}
}
\frac{1}{(m_1)^{i_1}(m_2)^{i_2}\cdots (m_k)^{i_k}}.\]
for symmetry with $t$-values.

\begin{fact}
Let $\{m\}^{n}$ denote the sequence 
$(\underbrace{m, m, \dots, m}_{n})$.
For a multi-index 
\[
\mathbf{i}=(a_1+1, 
\{1\}^{b_{1}-1},
a_2+1, \{1\}^{b_{2}-1}, \dots, a_k+1, 
\{1\}^{b_{k}-1}),
\]
with integers $k, a_j, b_j\ge1$, define its \emph{dual}  
\[
\mathbf{i}^\dagger=(b_k+1, 
\{1\}^{a_{k}-1}, 
b_{k-1}+1, \{1\}^{a_{k-1}-1},
 \dots, b_1+1, 
 \{1\}^{a_{1}-1}).
\]
\emph{Duality formula} for multiple zeta values claims  that 
$\z(\mathbf{i})=\z(\mathbf{i}^\dagger)$ for all indices such that the first argument is at least 2. 
However, there seems no duality formula for multiple $t$-values at the time of writing (December 2021). 
Historically, Drinfeld and Kontsevich found \emph{iterated integral expressions} for multiple zeta values and proved duality in 1990s. Afterward, Kaneko, Hoffman, Zagier and many other researchers developed the theory. There are also many applications of such sums to \eh{Euler sums} as discussed in the book \cite{valean}.
\end{fact}

\begin{ex}\hf
\begin{enumerate}
	\item 
the celebrated Euler-Goldbach theorem claims that 
\[
\z(2, 1)=\z(3).
\]
We can derive this relation 
from iterated integral expressions 
\[
\z(2, 1)=
\int_0^{1}\f{dx_3}{x_3}
\int_0^{x_3}\f{dx_2}{1-x_2}
\int_0^{x_2}\f{dx_1}{1-x_1}
\]
and 
\[
\z(3)=\int_0^{1}\f{dy_3}{y_3}
\int_0^{y_3}\f{dy_2}{y_2}
\int_0^{y_2}\f{dy_1}{1-y_1}
\]
with changing variables by $y_{j}=1-x_{4-j}$.
\item 
Observe that 
\[
t(2)=\ff{3}{4}\z(2),\q 
t(3)=\ff{7}{8}\z(3)
\]
and \cite[p.4]{hoffman} 
\[
t(2, 1)=-\f{1}{2}t(3)+t(2)\log2 \q (\ne t(3)).\]
\end{enumerate}
\end{ex}

\begin{thm}\label{thm4}
For $n\ge 1$, 
\[
\ff{\z(\{2\}^{n})}{2^{2n}}
=\ff{1}{(2n+1)!}\k{\ff{\pi}{2}}^{2n},
\q 
t(\{2\}^{n})
=\ff{1}{(2n)!}\k{\ff{\pi}{2}}^{2n}.
\]
\end{thm}
In fact, we can prove these by 
equating coefficients of 
$x^{2n}$ in 
\[
\dsum_{n=0}^{\mug}\ff{1}{(2n+1)!}
\k{\ff{\pix}{2}}^{2n}=
\ff{\sin \tff{\pix}{2}}{\tff{\pix}{2}}=
\prod_{n=1}^{\mug}\k{1-\ff{x^{2}}{(2n)^{2}}}
\]
and 
\[
\dsum_{n=0}^{\mug}\ff{1}{(2n)!}
\k{\ff{\pix}{2}}^{2n}=
\cos \ff{\pix}{2}=
\prod_{n=1}^{\mug}\k{1-\ff{x^{2}}{(2n-1)^{2}}}.
\]

However, we give a different proof here because it  suggests the application to evaluation of 
$\z(3, 2, \dots, 2)$ and $t(3, 2, \dots, 2)$ in the next subsection. For this purpose, we need a lemma.
\begin{lem}
\label{le3}
For $n\ge1$, $|x|\leqslant 1$, we have 
\[
\ff{\arcsin^{2n}x}{(2n)!}=
\sum_{
\substack{k>m_{1}>\cd >m_{n-1}>0\\}
}
\ff{2^{2k}}{\binom{2k}{k}}
\ff{1}{(2k)^{2}(2m_{1})^{2}\cd (2m_{n-1})^{2}}x^{2k},
\]
\[
\ff{\arcsin^{2n-1}x}{(2n-1)!}=
\sum_{
\substack{k>m_{1}>\cd >m_{n-1}\ge0\\}
}
\ff{\binom{2k}{k}}{2^{2k}}
\ff{1}{(2k+1)(2m_{1}+1)^{2}\cd (2m_{n-1}+1)^{2}}x^{2k+1}.
\]

\end{lem}

\begin{proof}
This is a rephrasing of 
J.M. Borwein--Chamberland {\cite[(1.1)-(1.4)]{boch}}. 
\end{proof}


\begin{proof}[Proof of Theorem \ref{thm4}]
Lemma \ref{le3} asserts that 
\[
\ff{\arcsin^{2n}x}{(2n)!}=
\sum_{
\substack{k>m_{1}>\cd >m_{n-1}>0\\}
}
\ff{2^{2k}}{\binom{2k}{k}}
\ff{1}{(2k)^{2}(2m_{1})^{2}\cd (2m_{n-1})^{2}}x^{2k}
\]
so that 
\[
W\k{\ff{\arcsin^{2n}x}{(2n)!}}
=
\f{\,\pi\,}{2}
\sum_{
\substack{k>m_{1}>\cd >m_{n-1}>0\\}
}
\ff{1}{(2k)^{2}(2m_{1})^{2}\cd (2m_{n-1})^{2}}x^{2k}.
\]
Let $x=1$. The left hand side becomes 
\[
\di\int_{0}^{1}{\ff{\arcsin^{2n}u}{(2n)!}}\,\ff{du}{\rt{1-u^{2}}}=
\ts{\ff{\arcsin^{2n+1}u}{(2n+1)!}}{1}{0}=
\ff{1}{(2n+1)!} \k{\f{\,\pi\,}{2}}^{2n+1}
\]
while the right hand side turns to be 
$\tff{\pi}{2}\tff{\z(\{2\}^{n})}{2^{2n}}$.
Hence we proved 
\[
\ff{\z(\{2\}^{n})}{2^{2n}}
=\ff{1}{(2n+1)!}\k{\ff{\pi}{2}}^{2n}.
\]
It is quite similar to show 
$t(\{2\}^{n})
=\tff{1}{(2n)!}\k{\tff{\pi}{2}}^{2n}$
using $\tf{\arcsin^{2n-1}x}{(2n-1)!}$ and the operator  $W$.
\end{proof}


\subsection{
$\z(3, 2, \dots, 2), t(3, 2, \dots, 2)$}

We just found $\z(2, \dots, 2)$ and $t(2, \dots, 2)$ above. 
A natural subsequence is to 
evaluate $\z(3, 2, \dots, 2)$ and $t(3, 2, \dots, 2)$ via  certain integration on arcsin.

\begin{defn}
For $n\ge1$, set 
\[I\left(n\right) = \int_0^1 \dfrac{\arcsin^nx}{x}\,dx.\]
\end{defn}


\begin{ex}
\label{arcsin_ints}
Observe the first several values.
\begin{equation}
\label{I1}
     I(1)= \dfrac{\pi}{2}\,\log\left(2\right)
\end{equation}
\begin{equation}
\label{I2}
     I(2)= \dfrac{\pi^2}{4}\,\log\left(2\right)-\dfrac{7}{8}\,\zeta(3)
    \end{equation}
\begin{equation}
     I(3)= \dfrac{\pi^3}{8}\,\log\left(2\right)-\dfrac{9\pi}{16}\,\zeta(3)
     \end{equation}
     \begin{equation}
     \label{I4}
     I(4)= \dfrac{\pi^4}{16}\,\log\left(2\right)-\dfrac{9\pi^2}{16}\,\zeta(3) + \dfrac{93}{32}\,\zeta(5)\end{equation}
     \begin{equation}
     I(5)= \dfrac{\pi^5}{32}\,\log\left(2\right) - \dfrac{15\pi^3}{32}\,\zeta(3) + \dfrac{225\pi}{64}\,\zeta(5)  \end{equation}
     \begin{equation}
     \label{I6}
I(6)=\dfrac{\pi^6}{64}\,\log\left(2\right) - \dfrac{45\pi^4}{128}\,\zeta(3) + \dfrac{675\pi^2}{128}\,\zeta(5) -\dfrac{5715}{256}\,\zeta(7) \end{equation}
\end{ex}

\begin{rmk}\hf
\begin{enumerate}
	\item Indeed, Wolfram alpha \cite{wolfram} returns the algebraic expressions (\ref{I1})-(\ref{I6}) for integrals 
\[
I(n)=\int_0^{\p/2} y^{n}\cot y\,dy
\]
while she outputs only numerical values for 
\[
\int_0^1 \dfrac{\arcsin^nx}{x}\,dx.
\]
In fact, there is a precise formula for $I(n)$ giving a 
$\qq$-linear combination of $\log2$ and single Riemann zeta values. 
For the sake of completeness, we discuss it here 
although we do not need it in the sequel.
Let $\eta$ denote the Dirichlet eta function, that is, 
$\eta\left(1\right)=\log2$ and 
$\eta\left(j\right) = \left(1-2^{1-j}\right)\zeta(j)$ $(j\ge2)$. 
Then, there holds
\begin{equation}
\label{I2N+1}
I\left(2n+1\right) = \dfrac{\left(2n+1\right)!}{2^{2n+1}}\sum_{j=0}^{n}\dfrac{\left(-1\right)^j\pi^{2n+1-2j}}{\left(2n+1-2j\right)!}\,\eta\left(2j+1\right), 
\end{equation}
\begin{equation}
\label{I2N}
I\left(2n\right) = \dfrac{\left(2n\right)!}{2^{2n}}
\left(\sum_{j=0}^{n-1}\dfrac{\left(-1\right)^j\pi^{2n-2j}}{\left(2n-2j\right)!}\,\eta\left(2j+1\right) + (-1)^n2(1-2^{-2n-1})\,\zeta\left(2n+1\right)\right).
\end{equation}
To see this, we remark that Buhler-Crandall \cite[p.280]{bucr} stated
\[
\displaystyle\int_{0}^{1/2}{x^{n}\cot(\pi x)}\,dx
\]
\[=\ff{n!}{2^{n}}
\sum_{
\substack{1\le k\le n\\
k \textrm{ odd}}
}
\ff{(-1)^{(k-1)/2}}{\pi^{k}}\ff{\eta(k)}{(n-k+1)!}
+\f12((-1)^{n}+1)\ff{4n!(1-2^{-n-1})}{(2\p)^{n+1}}\z(n+1).
\]
However, the sign $\f12((-1)^{n}+1)$ must be 
$\cos\tfrac{n\pi}{2}$ (for $n=2$, the coefficient of 
$\z(3)$ is negative; see (\ref{I2})).
To correct this, set 
\[
J(n):=
\displaystyle\int_{0}^{1/2}{x^{n}\cot(\pi x)}\,dx
\]
\[=
\ff{n!}{2^{n}}
\k{
\sum_{
\substack{1\le k\le n\\
k \textrm{\,odd}}
}
\ff{(-1)^{(k-1)/2}}{\pi^{k}}\ff{\eta(k)}{(n-k+1)!}
}
+
\cos\frac{n\pi}{2}
\ff{4n!(1-2^{-n-1})}{(2\p)^{n+1}}\z(n+1).
\]
Then, 
with $y=\sin(\pi x)$, we find 
\[
J(n)=\ff{1}{\pi^{n+1}}
\displaystyle\int_{0}^{1}{\ff{\arcsin^{n}(y)}{y}}\,dy
=\ff{1}{\p^{n+1}}I(n).
\]
Thus, $I(n)=\pi^{n+1}J(n)$.
Writing down the cases for the index even and odd with $k=2j+1$, we get 
(\ref{I2N+1}), (\ref{I2N})
 and hence justified (\ref{I1})-(\ref{I6}).
\item We can also view $I(n)$ as a \emph{log-sine  integral}:
\[
I(n)=
\int_0^1 \ff{\arcsin^n x}{x}dx=
\underbrace{
\left[
\log x (\arcsin^n x)
\right]_0^1}_{0}
-
n
\int_0^1 \log x \ff{\arcsin^{n-1} x}{\sqrt{1-x^2}}dx
\]
\[
=
-n
\int_0^{\pi/2} y^{n-1} \log(\sin y)dy.
\]
See J. M. Borwein-Broadhurst-Kamnitzer \cite{bobrka}  for relation of such integrals and central binomial series, for example. 
\end{enumerate}
\end{rmk}

%

\begin{thm}\label{thm5}
For $n\ge0$,
\[
\ff{\zeta(3, \{2\}^{n})}{2^{2n+3}}=
	\ff{2}{\pi}
\displaystyle\int_{0}^{1}
\ff{\arcsin^{2n+2}x}{(2n+2)!}
\ff{\text{\,arccos\,} x}{x}dx,
\]
\[
t(3, \{2\}^{n})=
\displaystyle\int_{0}^{1}{
\ff{\arcsin^{2n+1}x}{(2n+1)!}
\ff{\arccos x}{x}
}\,dx.
\]
\end{thm}
The proof is quite same to the one for Theorem  \ref{thm4}. Hence we omit it.

\begin{ex}
We already discussed the cases for $n=0$ as integrals 
giving $\tf{7}{8}\z(3)$ and $\tf{1}{8}\z(3)$.
For $n\ge1$, 
with $\arccos x=\tfrac{\pi}{2}-\arcsin x$, 
we see that 
\[
t(3, 2)=
\int_{0}^{1}
\frac{\arcsin^3x}{3!}
\ff{\arccos x}{x} dx =
\f{1}{3!}\k{\frac{\pi}{2}I(3)-I(4)}
=
\ff{1}{64}\k{
3\pi^{2}\z(3)-31\z(5)
},
\]
\[
\ff{\zeta(3,2)}{2^{5}}= 
\ff{2}{\p}
\f{1}{4!}\k{\frac{\pi}{2}I(4)-I(5)}
=
\ff{1}{64}
\k{\pi^{2}\z(3)-11\z(5)},
\]
\[
t(3, 2, 2)
=\f{1}{5!}\k{\frac{\pi}{2}I(5)-I(6)}
=
\f{1}{2048}
\k{2\pi^4\z(3)-60\pi^2\z(5)+381\z(7)}
\]
and so on.
\end{ex}

\begin{cor}
Let $\mathbf{Q}\left[\pi, \z(3), \z(5), \ldots, \z(2n+3)\right]_{2n+3}$
 denote the set of all elements of degree $2n+3$ in the rational polynomial ring in $\pi, \z(3), \z(5), \ldots, \z(2n+3)$ with grading $\deg \p=1$ and  $\deg \z(2j+1)=2j+1$. 
 Then 
\[
\ff{\z(3, \{2\}^{n})}{2^{2n+3}},  t(3, \{2\}^{n})
 \in 
\mathbf{Q}\left[\pi, \z(3), \z(5), \dots, \z(2n+3)\right]_{2n+3}.
\]
\end{cor}
\begin{rmk}
Not all multiple values satisfy such a property.
For example, as mentioned before, 
\[
t(2, 1)=\ff{\pi^{2}}{{8}}\log2-\ff{7}{16}\z(3)
\]
involves a rational multiple of $\pi^{2}\log2$.
\end{rmk}


\subsection{multiple mixed values}

We next discuss new kind of multiple sums.
%
%

\begin{defn}
For positive integers $i_1, \dots, i_k$ 
such that $i_{k}\ge2$, define two kinds of \emph{multiple mixed values}
\[
\mu\left(i_k, i_{k-1}, \ldots,i_1\right)
=\sum_{
\substack{n_k>n_{k-1}>\cdots >n_1\\
n_j \,\,\equiv \,\,j \,\,\textrm{(mod $2$)}
}}
\dfrac{1}{n_k^{i_k}n_{k-1}^{i_{k-1}} \ldots n_1^{i_1}},
\]
\[
\ol{\mu}\left(i_k, i_{k-1}, \ldots,i_1\right)
=\sum_{
\substack{n_k>n_{k-1}>\cdots >n_1\\
n_j \,\,\equiv \,\,j+1 \,\,\textrm{(mod $2$)}
}}
\dfrac{1}{n_k^{i_k}n_{k-1}^{i_{k-1}} \ldots n_1^{i_1}}.
\]
\end{defn}
\begin{ex}
\[
\mu(3)=t(3), \q \ol{\mu}(3)=\ff{\z(3)}{2^{3}}, 
\]
\[
\mu(2, 1)=\ff{1}{2}t(3), \q
\ol{\mu}(2, 1)=
t(3)-
t(2)\log2
\q \te{(as shown below)}.
\]
In particular, 
\[
\mu(3, 1)=
\sum_{
\substack{n_{2}>n_{1}\\
n_{1} \text{\, odd}\\
n_{2} \text{\, even}
}
}
\ff{1}{n_{2}^{3}n_{1}}
=0.16227...
\]
is known as \eh{Ramanujan constant} (often denoted by $G(1)$) \cite[p.255-257]{rama1}.
\end{ex}

\begin{rmk}
Since $(\arctanh x)'=\tf{1}{1-x^{2}}
=
\sum_{k=0}^{\mug}x^{2k}
$ for $|x|<1$, 
there exists the following iterated integral expression:
\[
\f{\text{arctanh}^{n+1} x}{(n+1)!}=
\int_0^x \f{dx_{n+1}}{{1-x_{n+1}^2}}
\cdots
\int_0^{x_{2}} \f{dx_1}{{1-x_1^2}}.
\]
Expanding each $\tff{1}{1-x_{j}^{2}}$ to geometric series, it holds that 
\begin{eqnarray}
\di\int_{0}^{1}{\ff{\arctanh^{n+1}x}{(n+1)!}}
\,
\ff{dx}{x}
=
\mu(2, \{1\}^{n}).
\label{eqmu}
\end{eqnarray}
\end{rmk}



We can naturally extend this little more. Let $n, r\ge0$ be integers. Notice that $(r+2, \{1\}^{n})^{\dagger}=(n+2, \{1\}^{r})$.


\begin{thm}\label{thm6}
\[
2^{n}\mu(r+2, \{1\}^{n})=
2^{r}\mu(n+2, \{1\}^{r}).
\]
\end{thm}

\begin{proof}
As a consequence of \kref{eqmu}, there exists an iterated integral 
expression 
\[
\int_0^{1} \f{dx_{n+r+2}}{{x_{n+r+2}}}
\cd
\int_0^{x_{n+3}} \f{dx_{n+2}}{{x_{n+2}}}
\int_0^{x_{n+2}} \f{dx_{n+1}}{{1-x_{n+1}^2}}
\cdots
\int_0^{x_{2}} \f{dx_1}{{1-x_1^2}}
=\mu(r+2, \{1\}^{n}).
%
\]
Now, introduce new variables $\{y_{j}\mid 1\le j\le n+r+2\}$ by 
\[
y_{j}=\ff{1-x_{j}}{1+x_{j}}.
\]
The map $x_{j}\mapsto y_{j}$ is an involution and 
order-reversing on $[0, 1]$.
Let us now see what happens to each of those integrals.
For $1\le j\le n+1$,
\[
\di\int_{0}^{x_{j+1}}{
\ff{dx_{j}}{1-x_{j}^{2}}
}=
\di\int_{1}^{y_{j+1}}{
\ff{1}{1-
\k{\tf{1-y_{j}}{1+y_{j}}}^{2}
}
\ff{-2}{(1+y_{j})^{2}}dy_{j}
}
=
\di\int_{y_{j+1}}^{1}{\ff{dy_{j}}{2y_{j}}}.
\]
For $n+2\le j\le n+r+2$, 
\[
\di\int_{0}^{x_{j+1}}{
\ff{dx_{j}}{x_{j}}
}
=
\di\int_{1}^{y_{j+1}}
\ff{1}{\tf{1-y_{j}}{1+y_{j}}}
\ff{-2}{(1+y_{j})^{2}}dy_{j}
=
\di\int_{y_{j+1}}^{1}{\ff{2dy_{j}}{1-y_{j}^{2}}}
\]
with $y_{n+r+3}=0$.
Altogether, 
\begin{align*}
	\mu(r+2, \{1\}^{n})&=
	\int_0^{1} \f{dx_{n+r+2}}{{x_{n+r+2}}}
\cd
\int_0^{x_{n+3}} \f{dx_{n+2}}{{x_{n+2}}}
\int_0^{x_{n+2}} \f{dx_{n+1}}{{1-x_{n+1}^2}}
\cdots
\int_0^{x_{2}} \f{dx_1}{{1-x_1^2}}
	\\&=
	\int_0^1 \f{dy_{1}}{2y_{1}}
\cdots
\int_0^{y_{n}} \f{dy_{n+1}}{2y_{n+1}}
\int_0^{y_{n+1}} \f{2dy_{n+2}}{{1-y_{n+2}^2}}
\cd
\int_0^{y_{n+r+1}} \f{2dy_{n+r+2}}{{1-y_{n+r+2}^2}}
\\
&=
\ff{2^{r+1}}{2^{n+1}}
\mu(n+2, \{1\}^{r}).
\end{align*}
\end{proof}

\begin{ex}
Theorem \ref{thm6} gives non-trivial relations among various multiple sums such as  
\[
2\mu(2, 1)=t(3), \q
4\mu(2, 1, 1)=t(4).
\]
Unfortunately, we cannot derive any relation on $G(1)=\mu(3, 1)$ since $(3, 1)^{\dagger}$ is itself.

Observe next that 
\[
\ol{\mu}(2, 1)=
\sum_{
\substack
{n>m\\
n:\te{odd}
\\m:\te{even}
}
}
\ff{1}{n^{2}m}=
\k{
\sum_{
\substack
{n>m\\
}
}
-
\sum_{
\substack
{n>m\\
n:\te{odd}
\\m:\te{odd}
}
}
-\sum_{
\substack
{n>m\\
n:\te{even}
\\m:\te{odd}
}
}
-
\sum_{
\substack
{n>m\\
n:\te{even}
\\m:\te{even}
}
}
}
\ff{1}{n^{2}m}
\]
\[=
\z(2, 1)-t(2, 1)-\mu(2, 1)-\ff{1}{8}\z(2, 1)
=
\z(3)-\k{
\ff{\pi^{2}}{8}\log2-\ff{7}{16}\z(3)}-
\ff{7}{16}\z(3)-\ff{1}{8}\z(3)
\]
\[
=\ff{7}{8}\z(3)-
\ff{\pi^{2}}{8}\log2.\]
Also, there holds the following simple relation which is not so obvious at a first glance:
\[
\mu(2, 2)+\ol{\mu}(2, 2)=2t(2, 2).
\]
Starting with 
\[
\dsum_{n=0}^{\mug}\ff{(-1)^{n}}{(2n+1)!}(\pix)^{2n}
=
\ff{\sin \px}{\px}=
\prod_{N=1}^{\mug}\k{1-\ff{x^{2}}{N^{2}}}
=
\prod_{m=1}^{\mug}\k{1-\ff{x^{2}}{(2m-1)^{2}}}
\prod_{n=1}^{\mug}\k{1-\ff{x^{2}}{(2n)^{2}}},
\]
equate the coefficients of $x^{4}$:
\[
\ff{1}{5!}\pi^{4}=
t(2, 2)+\ff{\z(2, 2)}{2^{4}}+\mu(2, 2)+\ol{\mu}(2, 2).
\]
Using 
$t(2, 2)=
\ff{1}{4!}\k{\ff{\pi}{2}}^{4}$,
$
\ff{\z(2, 2)}{2^{4}}=
\ff{1}{5!}\k{\ff{\pi}{2}}^{4}$,
we have 
\[
\mu(2, 2)+\ol{\mu}(2, 2)=\ff{1}{2^{3}4!}\pi^{4}=2t(2, 2).
\]
More generally, we have 
\[
\dsum_{n=0}^{\mug}\ff{(-1)^{m}}{(2m+1)!}(\pix)^{2m}
\dsum_{n=0}^{\mug}\ff{1}{(2n+1)!}(\pix)^{2n}
=
\ff{\sin \px}{\px}\ff{\sinh \px}{\px}
\]
\[
=
\prod_{N=1}^{\mug}\k{1-\ff{x^{4}}{N^{4}}}
=
\prod_{m=1}^{\mug}\k{1-\ff{x^{4}}{(2m-1)^{4}}}
\prod_{n=1}^{\mug}\k{1-\ff{x^{4}}{(2n)^{4}}}.
\]
Equate the coefficients of $x^{8}$:
\[
\k{\ff{1}{9!}-\ff{1}{3!7!}+\ff{1}{5!5!}-\ff{1}{7!3!}+
\ff{1}{9!}}\p^{8}=
t(4, 4)+\ff{\z(4, 4)}{2^{8}}
+\mu(4, 4)+\ol{\mu}(4, 4).
\]
With 
\[
t(4, 4)=\ff{\p^{8}}{2^{4}\cdot 8!}, \q
\ff{\z(4, 4)}{2^{8}}
=\ff{\p^{8}}{2^{3}\cdot 10!}, 
\]
($\z(\{4\}^{n})=\tff{2^{2n+1}\pi^{4n}}{(4n+2)!},
t(\{4\}^{n})=\tff{\pi^{4n}}{2^{2n}(4n)!}$ 
\cite[p.3]{hoffman}
), 
we finally have 
\[
\mu(4, 4)+\ol{\mu}(4, 4)=
\ff{\p^{8}}{138240}.
\]
\end{ex}


\subsection{conjecture}

We finish this article with some informal discussion without proofs. 
It is about powers of $\log x$ as certain operator; 
at the end, we make one conjecture.

For $f(x)\in \rr[[x]]$ such that $f(0)=0$, 
it is the technique to consider 
\[
\int_{0}^{y} \ff{f(x)}{x} dx
\]
to construct another series as we often encountered.

We can generalize this little more 
by changing the part ``$\int\tff{1}{x}$" with the iterated integral 
\[
\underbrace{\int 
\ff{dx}{x}
\int 
\ff{dx}{x}
\cd
\int 
\ff{dx}{x}}_{r+1}=
\int 
\ff{\log^{r}x}{r! x}.
\]

\begin{fact}[{\cite[p.1, 57-58]{valean}}]
For integers $n\ge1, r\ge0$, we have 
\[
\di\int_{0}^{1}{x^{n}\ff{\logrx}{r! x}}\,dx=
\ff{(-1)^{r}}{n^{r+1}}.
\]
\end{fact}
Now we obtain integral evaluations one after another:
\[
\di\int_{0}^{1}{{\athx}}
\ff{\logx}{x}
\,dx=
\di\int_{0}^{1}{
\dsum_{k=0}^{\mug}\ff{x^{2k+1}}{2k+1}
\ff{\logx}{x}
}dx
=
\dsum_{k=0}^{\mug}
\ff{1}{2k+1}
\di\int_{0}^{1}{x^{2k+1}\ff{\logx}{1!x}}\,dx
\]
\[
=
-
\dsum_{k=0}^{\mug}\ff{1}{(2k+1)^{3}}
=-\ff{7}{8}\z(3),
\]
\[
\di\int_{0}^{1}{{\athx}}
\ff{\log^{2}x}{2!x}
\,dx=
\ff{15}{16}\z(4),
\]
\[
\di\int_{0}^{1}{{\athx}}
\ff{\log^{3}x}{3!x}
\,dx=
-\ff{31}{32}\z(5),
\]
\[
\di\int_{0}^{1}{{\atx}}
\ff{\logx}{x}
\,dx=
-\ff{1}{32}\p^{3},
\]
\[
\di\int_{0}^{1}{{\atx}}
\ff{\log^{3}x}{3! x}
\,dx=-
\ff{5}{256}\p^{5},
\]
\[
\di\int_{0}^{1}{{\atx}}
\ff{\log^{5}x}{5! x}
\,dx=
-\ff{61}{1536}\p^{7}
\]
and more generally 
\[
\mu(r+2, \{1\}^{n})=
(-1)^{r}
\di\int_{0}^{1}{
\ff{\arctanh^{n+1}x}{(n+1)!}
\ff{\logrx}{r! x}
}\,dx.
\]
In particular, 
\[
G(1)=
\mu(3, 1)=
-\di\int_{0}^{1}{
\ff{\arctanh^{2}x}{2!}
\ff{\logx}{x}
}\,dx.
\]
Let us see what if we apply this idea back to $\arcsin$  integrals. Again, recall that 
\[
\asx
=\dsum_{k=0}^{\mug}\cbk
\ff{x^{2k+1}}{2k+1}.
\]
We now see another central binomial series
\[
-\di\int_{0}^{1}{\arcsin x\ff{\logx}{x}}\,dx=
\dsum_{k=0}^{\mug}
\cbk	
\ff{1}{{2k+1}}
\di\int_{0}^{1}{x^{2k+1}\ff{\logx}{x}}\,dx
=
\dsum_{k=0}^{\mug}
\cbk
\ff{1}{(2k+1)^{3}}.
\]
Value of this sum is known to be 
\[
\ff{1}{48}(\p^{3}+12\p\log^{2}2)\q 
\te{\cite[p.31]{rama2}}.
\] 
Similarly,
\[
\di\int_{0}^{1}{\arcsin x\ff{\log^{2}x}{2! x}}\,dx
=
\dsum_{k=0}^{\mug}
\cbk
\ff{1}{(2k+1)^{4}}.
\]
Wolfram alpha says that 
\[
\di\int_{0}^{1}{\arcsin x\ff{\log^{2}x}{2!x}}\,dx
=
\ff{1}{48}(6\p\z(3)+4\p\log^{3}2+\p^{3}\log2).
\]
This is only computer verification. 
Hence let us state it as a conjecture.
\begin{cj}
\[
\dsum_{k=0}^{\mug}
\cbk
\ff{1}{(2k+1)^{4}}
=
\ff{1}{48}(6\p\z(3)+4\p\log^{3}2+\p^{3}\log2).
\]
\end{cj}

\begin{rmk}
Recently, Ablinger (2015) \cite[p.17, 21]{ablinger} evaluated quite similar series 
\begin{align*}
	\dsum_{k=0}^{\mug}
\ff{\binom{2k}{k}}{4^{2k}}
\ff{1}{(2k+1)^{3}}
&=\ff{7}{216}\pi^{3},
	\\\dsum_{k=0}^{\mug}
\ff{\binom{2k}{k}}{4^{2k}}
\ff{1}{(2k+1)^{4}}
&=\ff{27\rt{3}}{32}
\dsum_{k=0}^{\mug}\ff{1}{(3k+1)^{4}}
+\ff{\pi}{12}\z(3)
-\ff{\p^{3}}{72\rt{3}}+\ff{27\rt{3}}{32}
\end{align*}
by method of iterated integrals, integration by parts, 
and generating functions. 
We expect that we can prove the conjecture with some similar idea.
\end{rmk}

\begin{center}
Acknowledgment.\\
The author would like to thank the anonymous referee for helpful comments to improve the manuscript. He also  
thanks Satomi Abe, Shoko Asami, Yuko Takada and 
Michihito Tobe for sincere support of his study.
\end{center}





%





%



\end{document}